\documentclass[letterpaper,12pt]{article}

\usepackage{physics}
\usepackage{tabularx} 
\usepackage{amsmath}  
\usepackage{graphicx} 
\usepackage[margin=1in,letterpaper]{geometry} 
\usepackage{cite} 
\usepackage[final]{hyperref} 
\hypersetup{
	colorlinks=true,       
	linkcolor=blue,        
	citecolor=blue,        
	filecolor=magenta,     
	urlcolor=blue         
}

\usepackage{subfig}

\usepackage{times}

\newenvironment{sciabstract}{%
\begin{quote}}
{\end{quote}}

\newenvironment{keywords}{%
\begin{quote}}
{\end{quote}}

\newcounter{lastnote}

\title{\bf Mathematical Construction of Interpolation and Extrapolation Function by Taylor Polynomials  } 

\author
{\bf Nijat Shukurov \\
\\
\normalsize{Department of Engineering Physics, Ankara University, Ankara, Turkey}\\
\\
\normalsize{E-mail:  nicatsukurov@gmail.com , nshukurov@ankara.edu.tr}
}

\date{}

\begin{document} 


\baselineskip24pt


\maketitle


\begin{sciabstract}
   \textbf{Abstract:} In this present paper, I propose a derivation of unified interpolation and extrapolation function that predicts new values inside and outside the given range by expanding direct Taylor series on the middle point of given data set. Mathematical construction of experimental model derived in general form. Trigonometric and Power functions adopted as test functions in the development of the vital aspects in numerical experiments. Experimental model was interpolated and extrapolated on data set that generated by test functions. The results of the numerical experiments which predicted by derived model compared with analytical values.
   
\end{sciabstract}

\begin{keywords}
	 \textbf{KEYWORDS:} Polynomial Interpolation, Extrapolation, Taylor Series, Expansion  

\end{keywords}

\newpage
\section{Introduction}

In scientific experiments or engineering applications, collected data are usually discrete in most cases and physical meaning is likely unpredictable. To estimate the outcomes and to understand the phenomena analytically controllable functions are desirable. In the mathematical field of numerical analysis those type of functions are called as interpolation and extrapolation functions.  Interpolation serves as the prediction tool within range of given discrete set, unlike interpolation, extrapolation functions designed to predict values out of the range of given data set. In this scientific paper, direct Taylor expansion is suggested as a instrument which estimates or approximates a new points inside and outside the range by known individual values. Taylor series is one of most beautiful analogies in mathematics, which make it possible to rewrite every smooth function as a infinite series of Taylor polynomials. As it is stated in \cite{1}, if the function $f(x)$ and its n derivatives are continuous within the range of special set which contains $x_{0}$ and $x_{a}$ points, then $f(x_{a})$ value can be detonated as.

\begin{align}
 f(x_{a})=f(x_{0})+(x_{a}-x_{0})f^{(1)}(x_{0})+ \frac{(x_{a}-x_{0})^{2}}{2!}f^{(2)}(x_{0}) \nonumber \\
  + ... + \frac{(x_{a}-x_{0})^{n}}{n!}f^{(n)}(x_{0}) + R_{n} 
\end{align} 

Reminder term $R_{n}$ defined as an integral form

\begin{align}
R_{n}=\frac{1}{n!}\int^{x}_{x_{0}} (x_{a}-x)^{n}f^{(n+1)}(x)dx   
\end{align} 

To understand the intuition behind the general analogy of proposed model it is useful to neglect the $R_{n}$ reminder term and just to focus on polynomial terms for now. Taylor's series is powerful mathematical tool which allows us to compute functions like trigonometric, logarithmic and exponential by approximation of these functions by Taylor polynomials. The aim of the proposed interpolation and extrapolation model is to approximate the Taylor polynomials of unknown function by feeding with known data set.

\section{Mathematical Construction of Proposed Model}

Consider N equidistant points in $x$ domain, which N is odd number. So that, $ \alpha=\frac{N-1}{2}$ and $x_{0}$ is the middle point of $\{ x_{k} \} \in [x_{-\alpha},x_{-\alpha+1},...,x_{\alpha+1},x_{\alpha}]$ set with step size of $\Delta x$. Also, there is N number of points in $y$ domain, which is the product of some unknown $y_{k}=f(x_{k})$ function and $\{x_{k},f(x_{k}) \}^{\alpha}_{k=-\alpha}$. In this case, $x_{k}=x_{0} + k \Delta x$ for $k \in [-\alpha , ..., \alpha]$  integer set. General formula of model is nothing but just Taylor expansion around $x_{0}$ point. 

\begin{align}
f^{model}_{0}(x)=\sum^{N}_{n=0} \frac{(x-x_{0})^{n} }{n!} \boldsymbol{ y^{(n)}_{0} }
\end{align}

To construct accurate Taylor expansion function for interpolation and extrapolation, it is desirable to generate $[y^{(0)}_{0},y^{(1)}_{0},...,y^{(N)}_{0}]$ set accurately. For this reason, $y_{k}=f(x_{0}+k\Delta x)$ unknown function was expanded via Taylor polynomials around the middle point $x_{0}$ for every element of $\{ y_{k} \}=\{ f(x_{0} + k \Delta x)\}$ set. As described in \cite{2}, to develop difference formulas it is necessary to use Taylor series.  

\begin{align}
y_{k}= y_{0}^{(0)} + (k \Delta x )^{1} y^{(1)}_{0}+ \frac{( k \Delta x )^{2}}{2!} y^{(2)}_{0} + ...
 + \frac{( k \Delta x )^{N-2}}{(N-2)!} y^{(N-2)}_{0}
 + \frac{( k \Delta x )^{N-1}}{(N-1)!} y^{(N-1)}_{0}  
\end{align}

Due to N points of data, same amount of Taylor expansions have produced matrix equation.

$$\begin{bmatrix}
   y_{-\alpha} \\
   y_{1-\alpha}\\
   \vdots\\
   y_{0}\\
  \vdots\\
   y_{\alpha-1}\\
   y_{\alpha}
  \end{bmatrix} = \begin{bmatrix} 
			1      & -\alpha \Delta x        & \cdots    &  \frac{[-\alpha\Delta x] ^{N/2}}{(N/2)! }       & \cdots & \frac{[-\alpha\Delta x] ^{N-2}}{(N-2)! }      & \frac{[-\alpha\Delta x] ^{N-1}}{(N-1)! }       \\ 
			1      & (1-\alpha) \Delta x        & \cdots & \frac{[(1-\alpha)\Delta x] ^{N/2}}{(N/2)! }      & \cdots & \frac{[(1-\alpha)\Delta x] ^{N-2}}{(N-2)! }       & \frac{[(1-\alpha)\Delta x] ^{N-1}}{(N-1)! }       \\ 
			\vdots & \vdots &        & \vdots  &        & \vdots & \vdots \\
			1      & 0      & \cdots & 0       & \cdots & 0      & 0      \\ 
			\vdots & \vdots &        & \vdots  &        & \vdots & \vdots \\			
			1      & (\alpha-1) \Delta x        & \cdots & \frac{[(\alpha-1)\Delta x] ^{N/2}}{(N/2)! }      & \cdots & \frac{[(\alpha-1)\Delta x] ^{N-2}}{(N-2)! }     & \frac{[(\alpha-1)\Delta x] ^{N-1}}{(N-1)! }       \\ 
			1      & \alpha \Delta x        & \cdots & \frac{[\alpha\Delta x] ^{N/2}}{(N/2)! }       & \cdots & \frac{[\alpha\Delta x] ^{N-2}}{(N-2)! }       & \frac{[\alpha\Delta x] ^{N-1}}{(N-1)! }      \\ 
							\end{bmatrix}   
							               \begin{bmatrix}
                                             y^{(0)}_{0} \\
  											 y^{(1)}_{0}\\
  												 \vdots\\
  											 y^{(N/2)}_{0}\\
  												 \vdots\\	 
   											 y^{(N-2)}_{0}\\
  											 y^{(N-1)}_{0}
  \end{bmatrix} $$

For simplicity matrix elements of matrix equation represented by capital letter as following.

\begin{align}
\boldsymbol{ Y_{N \times 1}} =\boldsymbol{ C_{ N\times N }} \boldsymbol{D_{N \times 1}}
\end{align}

Multiplying both sides by inverse of  $\boldsymbol{ C_{ N\times N }}$ matrix provides solution of  matrix equations, which is $ \boldsymbol{D_{N \times 1}} $ set that consist of all known order of derivations related with $y=f(x)$ function at $x=x_{0}$ point.

\begin{align}
\boldsymbol{C_{ N\times N }^{-1}} \boldsymbol{ C_{ N\times N }} \boldsymbol{D_{N \times 1}} \nonumber=&\boldsymbol{C_{ N\times N }^{-1}} \boldsymbol{ Y_{N \times 1}}  \\
\boldsymbol{D_{N \times 1}}=&\boldsymbol{C_{ N\times N }^{-1}}  \boldsymbol{Y_{N\times 1}}
\end{align}

Using that matrix elements, final formula of interpolation and extrapolation function has designed same as equation (3).
\begin{align}
f^{model}_{0}(x)=\sum^{N}_{n=0} \frac{(x-x_{0})^{n} }{n!} \boldsymbol{ D_{n} }
\end{align}

\section{Mathematical Expression of Numerical Experiment }
In this part, some numerical experiments have done with few points. First, let's introduce mathematical model for experimental reasons.

Consider $[x_{-2},x_{-1},x_{0},x_{1},x_{2}]$ and $[y_{-2},y_{-1},y_{0},y_{1},y_{2}]$ for $f(x_{i})=y_{i}$, and $\Delta x=x_{i+1}-x_{i}$ is a step distance of every point on $x$ domain.

Construct proposed interpolation and extrapolation function using direct Taylor expansion for given data around middle point $x_{0}$.
\begin{align}
f_{i}(x)=f^{(0)}(x_{0})+ (x-x_{0}) f^{(1)}(x_{0}) + \frac{(x-x_{0})^{2}}{2!}  f^{(2)}(x_{0}) \nonumber \\
 +  \frac{(x-x_{0})^{3}}{3!} f^{(3)}(x_{0})+\frac{(x-x_{0})^{4}}{4!} f^{(4)}(x_{0})
\end{align}

To achieve the interpolation and extrapolation function derivatives of unknown function should be known.
$$                                                      \begin{bmatrix}
														y_{-2}\\
														y_{-1}\\
														y_{0}\\
														y_{1}\\
														y_{2}
														\end{bmatrix}  =
            \begin{bmatrix} 
			1 & -2\Delta x & 2 \Delta x^{2}   &    \frac{-8 \Delta x^{3}}{3!}  &  \frac{16 \Delta x^{4}}{4!}   \\ 
			1 & -\Delta x & \frac{\Delta x^{2}}{2!}   &    \frac{- \Delta x^{3}}{3!}  &  \frac{\Delta x^{4}}{4!}   \\ 
			1 & 0 & 0 & 0 & 0 \\
			1 & \Delta x & \frac{\Delta x^{2}}{2!}   &    \frac{ \Delta x^{3}}{3!}  &  \frac{\Delta x^{4}}{4!}   \\ 
	1 & 2\Delta x & 2 \Delta x^{2}   &    \frac{8 \Delta x^{3}}{3!}  &  \frac{16 \Delta x^{4}}{4!}    
							\end{bmatrix}      
							               																													\begin{bmatrix}
														f^{(0)}(x_{0})\\
														f^{(1)}(x_{0})\\
														f^{(2)}(x_{0})\\
														f^{(3)}(x_{0})\\
														f^{(4)}(x_{0})
														\end{bmatrix}$$

\begin{align}
\boldsymbol{Y_{5 \times 1}} = \boldsymbol{C_{5 \times 5 }} \boldsymbol{D_{5 \times 1} } \\
\boldsymbol{D_{5 \times 1}} = \boldsymbol{C_{5 \times 5 }^{-1}} \boldsymbol{Y_{5 \times 1}}
\end{align}

The final form of unified interpolation and extrapolation functions is designed as following.

\begin{align}
f^{UIE}_{0}(x)=\boldsymbol{D_{1}}+ (x-x_{0}) \boldsymbol{D_{2}}+ \frac{(x-x_{0})^{2}}{2!} \boldsymbol{D_{3}}  \nonumber \\
 +  \frac{(x-x_{0})^{3}}{3!} \boldsymbol{D_{4}}+\frac{(x-x_{0})^{4}}{4!} \boldsymbol{D_{5}}
\end{align}

\newpage
\subsection{Numerical Experiments via Different Test Functions }

In this part using special set of x=[a,b] point y domain values was analytically calculated  respectively for every test functions to design experimental points for testing proposed interpolation and extrapolation function model. Using this data points model was set up and points within and outside the range of given set was generated and compared with analytically calculated points. These all numerical calculation have done by using \textbf{Python} programming language and the graphing of solution was done by using its \textbf{matplotlib} library. Algorithm of mathematical construction was implemented as a numerical Python module. Feed data was expected to be consumed by this module to design Taylor expansion function(proposed model). 

\subsubsection{$f(x)=ax^{3}+bx^{2}+cx+d$}
In this experiment to construct the experimental cubic function, $[a, b, c, d]$ coefficients was chosen randomly as $[3,2,1,4]$ respectively, so that function is $f(x)=3x^{3}+2x^{2}+x+4$.

\begin{table}[h]
\caption{} 
\centering 
\begin{tabular}{c rrrrr} 
\hline\hline 
Domain &\multicolumn{5}{c}{Feed Data} \\ [0.5ex]
\hline 
\textbf{x}  & -3.00   & -2.75   & -2.50  & -2.25  & -2.00  \\ 
\textbf{y}  & -62.00  & -46.01  & -32.88 & -22.30 & -14    \\[1ex] 
\hline 
\end{tabular}
\label{tab:table1}
\end{table}

$x$ domain and $y$ domain values (feed data) are represented at Table \ref{tab:table1} which is x data boundary has randomly chosen to generate 5 points with uniform $\Delta x=0.25$ step size and same amount of $y=f(x)$ values. In first step, to construct unified interpolation and extrapolation function it is desirable to find $\textbf{D}$ matrix (10), which generated by numerical Python model.  

$$\boldsymbol{D^{T}}=\begin{bmatrix}
														-32.87 & 47.25 & 41.00 & 18.00 & 0 \\
														
														\end{bmatrix}$$
                              
The final form of generated interpolation and extrapolation model (13) is nothing but just  Taylor expansion (12) of experimental cubic function at $x_{0}=-2.5$ point. 
\begin{align}
f^{Taylor}_{-2.5}(x)=-32.88+47.25(x+2.5)-20.50(x+2.5)^{2} \nonumber \\
+3.00(x+2.5)^{3}+O[(x+2.5)^{4}]
\end{align}
\begin{align}
f^{model}_{-2.5}(x)=&-32.88 + 47.25(x+2.5)-20.49(x+2.5)^{2}+3.00(x+2.5)^{3} \\
f^{model}_{-2.5}(x) \approx &f_{-2.5}^{Taylor}(x) 
\end{align}

To make it understandable visually, proposed interpolation and extrapolation model function values scattered by green triangles within [-5,5] range  which contains test points itself. Also, feed data and analytical values of cubic function has been scattered and plotted with blue circles and yellow squares respectively.

\begin{figure}[h]
\centering
\includegraphics[scale=0.7]{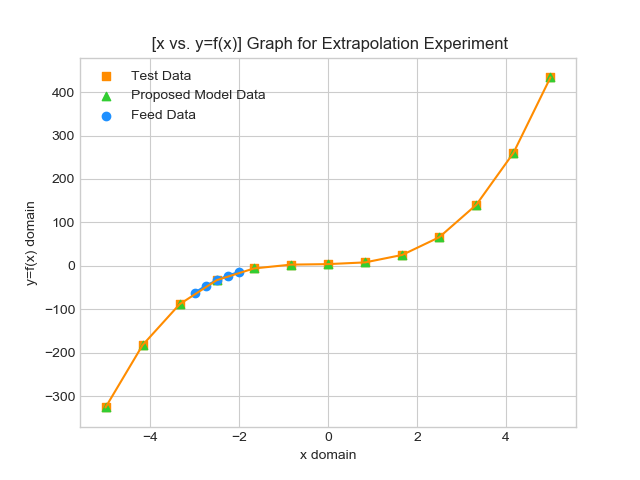}
\caption{Cubic Test}
\label{fig:figure1}
\end{figure}

Due to numerical power of Taylor expansion, there is no almost value differences between model and test cubic function in given range. To make difference visible or to catch the error it is convenient to use very large values which represented in Table \ref{tab:table2}.

\begin{table}[h!]
  \begin{center}
    \caption{}
    \label{tab:table2}
    \begin{tabular}{c rrrr}
      $x_{i}$ &  Cubic function: $f(x_{i})$   &   Proposed model: $f^{model}_{-2.5}(x_{i})$ &  $Err.=f(x_{i})-f^{model}_{-2.5}(x_{i})$ \\ 
      \hline
      999 & $2.99300600206\times10^{9}$ & $2.993006002\times10^{9}$ &    -0.0570 \\ 
      9999 & $2.999300060002\times10^{12}$ & $2.99930006057\times10^{12}$ & -568.6806\\ 
    \end{tabular}
  \end{center}
\end{table}

\subsubsection{$f(x)=ax^{4}+bx^{3}+cx^{2}+dx+f$}
Same as part before, to construct the experimental quadratic function, coefficients has choosen randomly, so the final form is $f(x)=5x^{4}+3x^{3}+1x^{2}+4x+2$.
\begin{table}[h]
\caption{} 
\centering 
\begin{tabular}{c rrrrrrr} 
\hline\hline 
Domain &\multicolumn{5}{c}{Points and Values} \\ [0.5ex]
\hline 
\textbf{x}  & 3.00 & 3.25 & 3.50 & 3.75 & 4.00 \\ 
\textbf{y}  & 509.00  & 686.38  & 907.19 & 1178.04 & 1506.00 \\[1ex] 
\hline 
\end{tabular}
\label{tab:table3}
\end{table}

This experiment details follows the same steps with cubic function experiment. Following Table \ref{tab:table3} represents $x$ domain range [-3,4] which contains 5 equadistant points with $\Delta x=0.25$ and generated test values. Interpolation and extrapolation function (16) constructed using $\textbf{D}$ matrix, and Taylor expansion of test quadratic function (15) are represented respectively below.
$$\boldsymbol{D^{T}}=\begin{bmatrix}
														907.19 & 978.75 & 800.00 & 438.00 & 120.00 \\
														
														\end{bmatrix}$$
                              
\begin{align}
f^{Taylor}_{3.5}(x)= & 907.19+978.75(x-3.5) +400.00(x-3.5)^{2}+73.00(x-3.5)^{3}+5.00(x-3.5)^{4}\nonumber \\
          & \quad +O[(x+2.5)^{5}]  \\
f^{model}_{3.5}(x)=&907.19 + 978.74(x-3.5) + 400.00(x-3.5)^{2}+73.00(x-3.5)^{3} \nonumber\\
          &+5.00(x-3.5)^{4}\\
f^{model}_{3.5}(x) \approx &f^{Taylor}_{3.5}(x) 
\end{align}

\begin{figure}[h]
\centering
\includegraphics[scale=0.7]{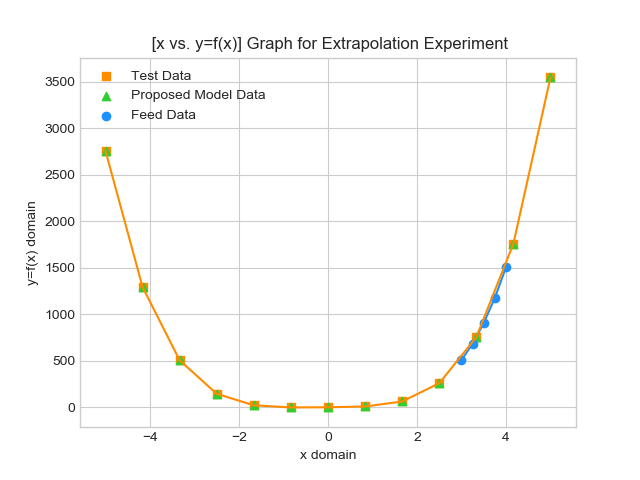}
\caption{Numerical Experiment by Quadratic Function values.}
\label{fig:figure2}
\end{figure}

Using all those relation following functions was plotted at Figure \ref{fig:figure2}.

Similarly cubic function experiment, it is severe to observe error visually at quadratic function graphic too. For that reason it is desirable to use very large values which is represented in Table \ref{tab:table4}.

\begin{table}[h!]
  \begin{center}
    \caption{}
    \label{tab:table4}
    \begin{tabular}{c rrrr}
      $x_{i}$ &  Quadratic function: $f(x_{i})$   &   Proposed model: $f^{model}_{3.5}(x_{i})$ &  $Err.=f(x_{i})-f^{model}_{3.5}(x_{i})$ \\ 
      \hline
      999 & $4.983021991001\times 10^{12}$ & $4.983021991\times10^{12}$ &    -0.0019 \\ 
      9999 & $4.9983002199910001\times 10^{16}$ & $4.99830021999\times10^{16}$ & 10001.0000\\ 
    \end{tabular}
  \end{center}
\end{table}

\subsubsection{$f(x)=sin(x)$}
\begin{figure}[h]
\centering
  \subfloat[]{\includegraphics[width=0.5\textwidth]{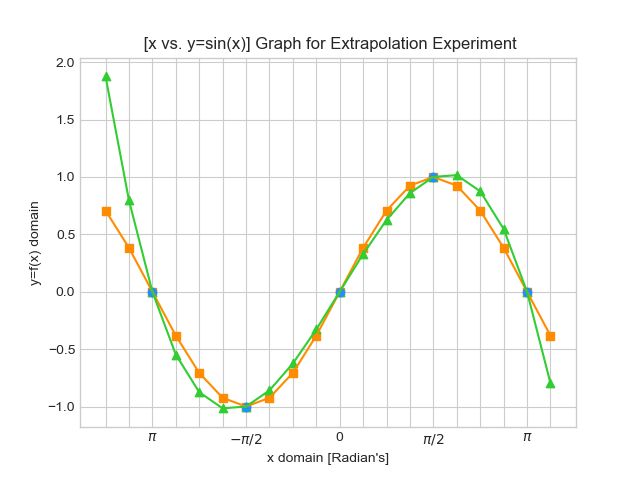}\label{fig:figure3a}}
  \hfill
  \subfloat[]{\includegraphics[width=0.5\textwidth]{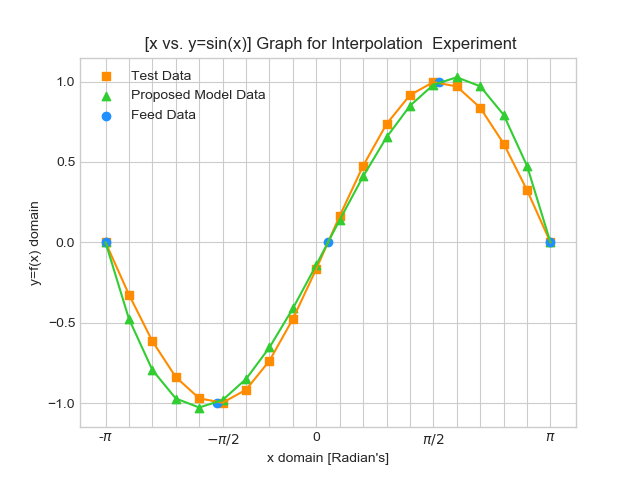}\label{fig:figure3b}}
  \caption{{These graphics are generated by 5 point feed data model. a) Extrapolation, b) Interpolation }}
\end{figure}

In this section, subject undergoes to smooth trigonometric function test in two different paths, interpolations and extrapolation. Firstly, same as the sections before, feed values have generated by randomly chosen $[-\pi,\pi]$ range on $x$ domain which contains 5 points with same step-size. These data have fed by mathematical model to achieve Taylor expansion. All following steps are similar with power function experiments. Outcomes was clearly shown in Figure \ref{fig:figure3a}, Figure \ref{fig:figure3b} and Table \ref{tab:table5}.

\begin{figure}[h]
\centering
  \subfloat[]{\includegraphics[width=0.5\textwidth]{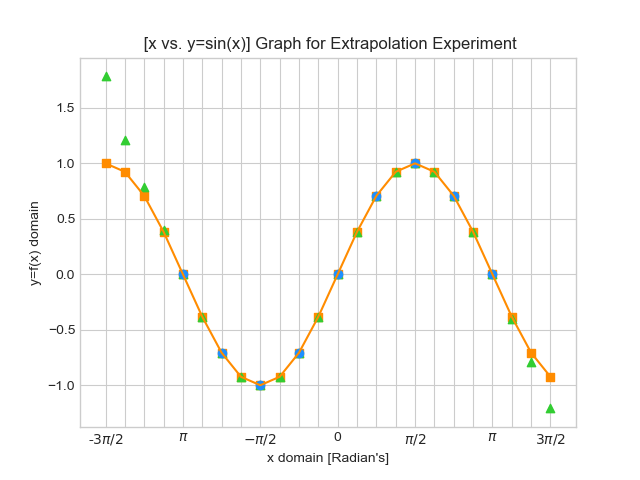}\label{fig:figure4a}}
  \hfill
  \subfloat[]{\includegraphics[width=0.5\textwidth]{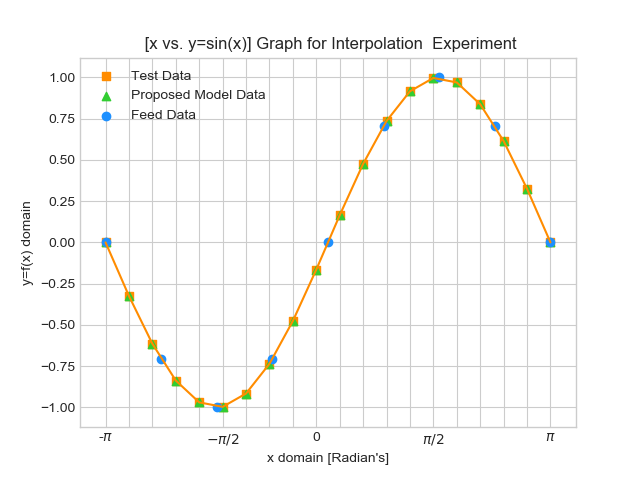}\label{fig:figure4b}}
  \caption{These graphics are generated by 9 point feed data model. a) Extrapolation, b) Interpolation }
\end{figure}

It is visible that this model is not pretty successful with smooth trigonometric functions in extrapolation and interpolation by using small number of feed data. For that reason using 9 points of data has fed by model then tested and results are little satisfying Figure \ref{fig:figure4a} and \ref{fig:figure4b}. Still proposed model is not successful for extrapolation to predict values outside the given range Figure \ref{fig:figure4a}. Interpolation \ref{fig:figure4b} feature of this model is quite accomplishing, there is still less accuracy in predicting values inside the region as shown in Table \ref{tab:table5}. 

\begin{table}[h!]
  \begin{center}
    \caption{5 points feed data}
    \label{tab:table5}
    \begin{tabular}{c rrrr}
      $x_{i}$ &  Sinusoidal function: $f(x_{i})$   &   Proposed model: $f^{model}_{0}(x_{i})$ &  $Err.=f(x_{i})-f^{model}_{0}(x_{i})$ \\ 
      \hline
      $2 \pi / 3 $ &  0.8660 & 0.9876 &    -0.1216 \\ 
      $5 \pi / 6 $ & 0.4999 & 0.6790 & -0.1790\\ %
    \end{tabular}
  \end{center}
\end{table}

\begin{table}[h!]
  \begin{center}
    \caption{9 points feed data}
    \label{tab:table6}
    \begin{tabular}{c rrrr}
      $x_{i}$ &  Sinusoidal function: $f(x_{i})$   &   Proposed model: $f^{model}_{0}(x_{i})$ &  $Err.=f(x_{i})-f^{model}_{0}(x_{i})$ \\ 
      \hline
      $2 \pi / 3 $ &  0.8660 & 0.8658 &    0.0002 \\ 
      $5 \pi / 6 $ & 0.4999 & 0.5006 & -0.0006\\ %
    \end{tabular}
  \end{center}
\end{table}

\section{Discussions}
Based on the results listed in Table \ref{tab:table1} and Table \ref{tab:table2}, interpolation and extrapolation function by Taylor polynomials provide accuracy in power functions, but there is dramatic increase in errors at very distant points because of Taylor approximation characteristics. Unlike power functions, this numerical module is not very precise at smooth trigonometric functions which is unreliable with less number of feed data as shown at Figure \ref{fig:figure3a} and Figure \ref{fig:figure3b}. In case of increased amount of feed data improved accuracy was recorded at interpolating process. In contrast to interpolation, there is still considerable more error in extrapolating Figure \ref{fig:figure4a}, and extrapolation feature of this model is not reliable for sinusoidal type data. But, interpolation specialty is well grounded due to result at Table \ref{tab:table6}. As demonstrated in graphic on Figure \ref{fig:figure4b} derived model is utilizable at interpolating by numerous sinusoidal feed data. To sum up, interpolation and extrapolation function by Taylor polynomials is excellent tool for constructing new points within the known sequence and estimation of a points based on extending a given set.

\section{Conclusion} 
Mathematical construction of interpolation and extrapolation function by Taylor polynomials was derived in general form. This model developed for experimental problems, and algorithm implemented in Python programming language as numerical module. The accuracy of proposed model was tested in three different numerical experiments, and results were compared with analytical values.

\end{document}